\documentclass[12pt]{article}

\usepackage{enumerate}
\usepackage{amsmath}
\usepackage{amssymb,latexsym}
\usepackage{amsthm}

\usepackage{anysize}
\usepackage{graphicx}
\usepackage{url}
\usepackage{color}
\usepackage[english]{babel}
\usepackage{latexsym}
\usepackage{amssymb,amsthm,amsmath}
\usepackage{graphicx,color}
\usepackage{verbatim,setspace}
\newtheorem{theorem}{Theorem}[section]

\newtheorem{lemma}[theorem]{Lemma}

\newtheorem{remark}[theorem]{Remark}
\newcommand{\Aut}{{\rm Aut}}
\newcommand{\Paut}{{\rm Paut}}

\newcommand{\diag}{{\rm diag}}
\newcommand{\ord}{{\rm ord}}

\newcommand{\QQ}{{\bf Q}}

\newcommand{\Z}{{\bf Z}}

\newcommand{\SSS}[1]{\mbox{${\bf S}_{#1}$}}

\newcommand{\Proof}{{\it Proof}.\ }

\begin{document}

\title{The First and Second Most Symmetric Nonsingular Cubic Surfaces
}
\author{H. Kaneta, S. Marcugini and F.Pambianco}

\maketitle

\begin{abstract}
The first and second most symmetric nonsingular cubic surfaces are $x^3+y^3+z^3+t^3=0$ and $x^2y+y^2z+z^2t+t^2x=0$, respectively.

\end{abstract}

\maketitle
\section{Introduction}

Throughout this paper $k$ stands for an algebraically closed field of characteristic zero. Let $C_d$ be a nonsingular curve of degree $d\geq 4$
in the two-dimensional projective space. As is well known, its automorphism group $\Aut(C_d)$ is a finite subgroup of $PGL_3(k)$ such that
$|\Aut(C_d)|\leq 42d(d-3)$. In fact the maximum value $\beta_d$ of $|\Aut(C_d)|$ is equal to $168$, $360$ or $6d^2$, and attained by the Klein
quartic curve, the Wiman sextic curve or the Fermat curve according as $d=4$, $d=6$ or $d\not\in \{4,6\}$\cite{pam}. Moreover the Klein quartic
and the Wiman sextic give rise to highly symmetric MDS codes \cite{kmp}.

Let $S_d$ be a nonsingular surface of degree $d\geq 3$ in the three-dimensional projective space. Its automorphism group
$\Aut(S_d)$ is a finite subgroup of $PGL_4(k)$ if $d\not=4$, can be an infinite group if $d=4$ \cite{mat}.
Highly symmetric surfaces may be both of theoretical and of practical interest.
Hosoh has completed the classification of
the automorphism groups of nonsingular cubic surfaces \cite{hos}. According to his classification  the maximum and the second maximum of $|\Aut(S_3)|$
is attained by the semidirect product $(\Z_3)^3\times_s \SSS{4}$ and $\SSS{5}$, respectively. It may happen that $|\Paut(S_d)|$, the order of the projective
automorphism group of $S_d$, has the finite upper bound for every $d\geq 3$. In this paper we will specify $S_3$ such that $\Aut(S_3)$ is isomorphic to
$(\Z_3)^3\times_s \SSS{4}$ or $\SSS{5}$.

$M_{m,n}(k)$ stands for the set of all $m\times n$ matrices of entries in $k$.
By definition $M_n(k)=M_{n,n}(k)$, $GL_n(k)=\{A=[a_{ij}]\in M_n(k)\ :\ \det A\not=0\}$, and $PGL_n(k)=GL_n(k)/( E_n)$,
where $(E_n)$ is the subgroup $\{\lambda E_n\ : \ \lambda\in k^*\}$ ($E_n$ is the unit matrix in $GL_n(k)$).
The coset $A(E_n)$ containing an $A\in GL_n(k)$ will be  denoted $(A)$. We denote by $k[x]$ the $k$-algebra of polynomials in $x=[x_1,...,x_n]$ over $k$.
For an $A\in GL_n(k)$ and $f\in k[x]$ we define a polynomial $f_A\in k[x]$ to be $f_A(x)=f(\sum \alpha_{1j}x_j,...,\sum \alpha_{nj}x_j)$, where
$A^{-1}=[\alpha_{ij}]$. As  is well known, the map $T_A:k[x]\rightarrow k[x]$ assigning $f_A$ to $f$ is a $k$-algebra isomorphism of $k[x]$ such that
$T_AT_B=T_{AB}$, that is, $(f_B)_A=f_{AB}$. $\SSS{n}$ denotes the symmetric group, the group of all permutations of $n$ elements.

A homogeneous polynomial $f$ of degree $d\geq 1$ defines a projective algebraic set $V(f)=\{(a)\in P^{n-1}\ : \ f(a)=0\}$ of an $(n-1)-$dimensional
projective space $P^{n-1}$ over $k$. $V(f)$ is called a hypersurface of degree $d$. If $n=4$, $d=3$, and $V(f)$ is nonsingular, the automorphism group
 $\Aut(V(f))$ of the hypersurface $V(f)$ is a finite subgroup of $PGL_4(k)$ consisting of $(A)\in PGL_4(k)$ such that $f_{A}\sim f$ \cite{mat}.
Let $a=[a_1,...,a_n]\in k^n$, $(a)\in V(f)$, and $A\in GL_n(k)$. Then $(a)$ is a singular point of $V(f)$ if $f_{x_i}(a)=0$ for all $i$. If
$(a)$ is a nonsingular point of $V(f)$,  $V(\sum_{i=1}^n\gamma_ix_i)$ is the tangent plane to $V(f)$ at $(a)$, where $\gamma_i=f_{x_i}(a)$. Clearly
$(A):V(f)\rightarrow V(f_{A})$ is a bijection, and if $b=Aa$ with $(a)\in V(f)$, then $(f_A)_{x_j}(b)=\sum_{i=1}^n f_{x_i}(a)\alpha_{ij}$, where
$A^{-1}=[\alpha_{ij}]$. Consequently $(b)=(A)(a)$ is a nonsingular point of $V(f_A)$ if and only if $(a)$ is a nonsingular point of $V(f)$, and the tangent
plane of $V(f_A)$ at $(b)$ coincides with $(A)V(\sum_{i=1}^n \gamma_ix_i)$. In particular if $(a)$ is a nonsingular point of $V(f)$, $f_{A}\sim f$ and
$(A)(a)=(a)$, then $[f_{x_1}(a),...,f_{x_n}(a)]A\sim [f_{x_1}(a),...,f_{x_n}(a)]$.

In $\S 1$ it will be shown that the Fermat surface $S=V(x^3+y^3+z^3+t^3)$ is the unique cubic nonsingular surface, up to projective equivalence, such that
$\Aut(S)$ is isomorphic to  $(\Z_3)^3\times_s \SSS{4}$.
In $\S 2$ it will be shown that a surface $S'=V(x^2t+y^2z+z^2t+t^2x)$ is the unique
cubic nonsingular surface such that $\Aut(S')$ is isomorphic to $\SSS{5}$.
\section{${\Z_3}^3\times_{s}\SSS{4}$-invariant nonsingular cubic surfaces}
Let $\omega\in k^*$ be of order three. Any subgroup of $PGL_4(k)$ isomorphic to $\Z_3$ is conjugate to $\langle (\diag[\omega,1,1,1])\rangle$ or
$\langle (\diag[\omega,\omega^2,1,1])\rangle$. It can be verified easily that $X=[x_{ij}]\in GL_4(k)$ satisfies
$\diag[\omega,1,1,1]X\sim X\diag[\omega,1,1,1]$ (resp. $\diag[\omega,\omega^2,1,1]X\sim X\diag[\omega,\omega^2,1,1]$) if and only if
$x_{1i}=x_{i1}=0$ for all $i\in [2,4]$ (resp. $x_{1i}=x_{i1}=x_{2j}=x_{j2}=0$ for all $i\in [2,4]$ and all $j\in [3,4]$). Consequently
a subgroup of $PGL_4(k)$ isomorphic to $(\Z_3)^2$ is conjugate to $\langle (A),(B)\rangle$, where $\ord(A)=\ord(B)=3$ and $(B^j)\in \langle (A)\rangle$
if and only if $j\in 3\Z$. We may assume that $A=\diag[\omega,\omega^2,1,1]$ or $A=\diag[\omega,1,1,1]$ and that $B=\diag[b_1,b_2,b_3,b_4]$.
Assume first that $A=\diag[\omega,\omega^2,1,1]$.
If $|\{b_1,b_2,b_3,b_4\}|=3$, then we may assume that $B$ is equal to one of $S-\{A\}$, where
\begin{eqnarray*}
 S&=&\{\diag[1,1,\omega,\omega^2], \diag[1,\omega,1,\omega^2], \diag[1,\omega,\omega^2,1], \diag[\omega,1,1,\omega^2], \\
    && \diag[\omega,1,\omega^2,1], \diag[\omega,\omega^2,1,1] \},
\end{eqnarray*}
Unless $B=\diag[1,1,\omega,\omega^2]$, there exist integers
$i$ and $j$ such that $(A^iB^j)$ is equal to $(\diag[1,1,\omega,1])$ or $(\diag[1,1,1,\omega]).$  So $\langle (A),(B)\rangle$ is conjugate to
$$\langle (\diag[\omega,\omega^2,1,1]),(\diag[1,1,\omega,\omega^2])\rangle \;\; or \;\; \langle (\diag[\omega,\omega^2,1,1]),(\diag[1,1,\omega,1])\rangle,$$
provided $|\{b_1,b_2,b_3,b_4\}|=3$.
Clearly $\langle (A),(B)\rangle$ is conjugate to $$
\langle (\diag[\omega,\omega^2,1,1]),(\diag[1,1,\omega,1)\rangle \;\; or \;\; \langle (\diag[\omega,1,1,1]),(\diag[1,1,\omega,1])\rangle,$$
provided $|\{b_1,b_2,b_3,b_4\}|=2$. Assume secondly that $A=\diag[\omega,1,1,1]$.  According as $|\{b_1,b_2,b_3,b_4\}|$ is two or three, we may
assume $B$ is equal to one of $T-\{A\}$ or $S$, where
\begin{eqnarray*}
T&=&\{\diag[\omega,1,1,1], \diag[1,\omega,1,1], \diag[1,1,\omega,1], \diag[1,1,1,\omega]\}.
\end{eqnarray*}
Thus  a subgroup of $PGL_4(k)$ isomorphic to $(\Z_3)^2$ is conjugate to one of
 $$\langle (\diag[\omega,\omega^2,1,1]),(\diag[1,1,\omega,\omega^2)]\rangle,\ $$
$$\langle (\diag[\omega,\omega^2,1,1]),(\diag[1,1,\omega,1])\rangle,\ $$
$$\langle(\diag[\omega,1,1,]),(\diag[1,\omega,1,1])\rangle.$$
Even though the following lemma may be known, we shall give a proof for the sake of completeness.
\begin{lemma}
A subgroup of $PGL_4(k)$ isomorphic to $(\Z_3)^3$ is conjugate to
\[
G_{27}=\langle(\diag[\omega,1,1,1]),(\diag[1,\omega,1,1]),(\diag[1,1,\omega,1])\rangle.
\]
\end{lemma}
\Proof
Let $G=\langle(A),(B),(C)\rangle$ be a subgroup of $PGL_4(k)$ isomorphic to $(\Z_3)^3$. Since $\langle(A),(B)\rangle$ is isomorphic to $(\Z_3)^2$,
we may assume one of the following three cases:
1) $A=\diag[\omega,\omega^2,1,1],\ B=\diag[1,1,\omega,\omega^2].$
2) $A=\diag[\omega,\omega^2,1,1],\ B=\diag[1,1,\omega,1].$
3) $A=\diag[\omega,1,1,1],\ B=\diag[1,\omega,1,1].$
In any case we may assume $C=\diag[c_1,c_2,c_3,c_4]\not=E_4$ with $c_i^3=1$ but $(C)\not\in \langle(A),(B)\rangle$.
Let $m=|\{c_1,c_2,c_3,c_4\}|$.
Assume the case 1).
If $m=2$, then $G=G_{27}$. If $m=3$, then we may assume $C$ is equal to one of $S-\{A,B\}$
so that $G$ contains one of $T$, hence $G=G_{27}$.
Assume the case 2).
If $m=2$, then we may assume $C$ is one of $T-\{B\}$, hence $G=G_{27}$. If $m=3$, then we may assume $C$ is one of
$S-\{ \diag[1,\omega,1,\omega^2], \diag[\omega,1,1,\omega^2], A\}$, hence $G=G_{27}$. Note that
$(\diag[1,\omega,1,\omega^2]),\ (\diag[\omega,1,1,\omega^2])\in \langle(A),(B)\rangle$.
Assume the case 3). If $m=2$, then
we may assume $C$ is one of $T-\{A,B\}$, hence $G=G_{27}$. If $m=3$, then we may assume $C$ is one of $S-\{\diag[\omega,\omega^2,1,1]\}$, hence $G=G_{27}$. \\ 

The canonical group representation $\hat{}:\SSS{4}\rightarrow GL_4(k)$ of $\SSS{4}$ is the one such that $\hat{\sigma}x=y$ with $y_i=x_{\sigma^{-1}(i)}$ for
any column vector $x\in k^4$.
Clearly this representation is an isomorphism. Let $r\geq 2$ be an integer, $\delta\in k^*$ be of order $r$,
\[
 D(r)=\{\diag[\lambda_1,\lambda_2,\lambda_3,\lambda_4]\in GL_4(k)\ :\ \lambda_1^r=\lambda_2^r=\lambda_3^r=\lambda_4^r\}.
\]
The factor group $(D(r))=D(r)/k^*E_4$ is isomorphic to $\{\diag[\delta^i,\delta^j,\delta^\ell,1]\ :\ i,j,\ell\in [0,r-1]\}$ which is isomorphic to $(\Z_r)^3$.
A map $\varphi_\sigma:(D(r))\rightarrow (D(r))$ defined by $\varphi_\sigma((A))=(\hat{\sigma}A{\hat{\sigma}}^{-1})$ is a group automorphism such that
$\varphi_{\sigma\tau}=\varphi_{\sigma}\circ\varphi_{\tau}$. The factor group $D(r) \hat {{\bf S}}_4/k^*E_4$ is isomorphic to the semidirect product
$(D(r))\times_s \SSS{4}$ such that $((A),\sigma)((A'),\sigma')=((A)\varphi_\sigma(A'),\sigma\sigma')$. Clearly $(D(r))\times_s \SSS{4}$ is isomorphic to
$(\Z_r)^3\times_s \SSS{4}$. We note that $\hat{\sigma}\diag[a_1,a_2,a_3,a_4]{\hat{\sigma}}^{-1}=\diag[b_1,b_2,b_3,b_4]$, where $b_i=a_{\sigma^{-1}(i)}$.
We may skip the proof of the following lemma \cite{shi}.
\begin{lemma}If $r\geq 3$, then the projective automorphism group of the surface $V(x^r+y^r+z^r+t^r)$ is $D(r)\hat{{\bf S}}_4/k^*E_4$, which can be
identified with $(D(r))\times_s \SSS{4}$.
\end{lemma}
\begin{lemma} Any $(\Z_3)^3$-invariant nonsingular cubic surface is projectively equivalent to $V(x^3+y^3+z^3+t^3)$. In particular
any $(\Z_3)^3\times_s \SSS{4}$-invariant nonsingular cubic surface is projectively equivalent to $V(x^3+y^3+z^3+t^3)$, and its automorphism group is
conjugate to $(D(3))\times_s \SSS{4}$.
\end{lemma}
\Proof By Lemma 1.1 it suffices to show that any $G_{27}$-invariant nonsingular cubic surface is $V(ax^3+by^3+cz^3+dt^3)$, where $a$, $b$, $c$ and $d$ are
nonzero constants. Let $V(f)$ be a $G_{27}$-invariant nonsingular cubic surface, where the homogeneous polynomial $f(x,y,z,t)$ of degree three has the form
\begin{eqnarray*}
&&a_1x^3+a_2y^3+a_3z^3+a_4t^3\\
&&+x^2(b_{12}y+b_{13}z+b_{14}t)+y^2(b_{21}x+b_{23}z+b_{24}t)+z^2(b_{31}x+b_{32}y +b_{34}t)\\ &&+t^2(b_{41}x+b_{42}y+b_{43}z)+c_1yzt+c_2xzt+c_3xyt+c_4xyz.
\end{eqnarray*}
Let
\[
 A_1=\diag[\omega,1,1,1],\ A_2=\diag[1,\omega,1,1],\ A_3=\diag[1,1,\omega,1].
\]
It is evident that $\ord(A_i)=3$.
Since $G_{27}$ contains $(A_i)$ ($i\in [1,4]$), $f_{A_i^{-1}}$ is equal to one of $\{f,\ \omega f,\ \omega^2f\}$. If $f_{A_1^{-1}}$ is equal to
$\omega f$ or $\omega^2f$, then $f$ is divisible by $x$, hence $V(f)$ is singular. Assume $f_{A_1^{-1}}=f$. Then
\[
 f(x,y,z,t)=a_1x^3+a_2y^3+a_3z^3+a_4t^3+y^2(b_{23}z+b_{24}t)+z^2(b_{32}y+b_{34}t)+t^2(b_{42}y+b_{43}z)+c_1yzt.
\]

 Unless $f_{A_2^{-1}}=f$, $V(f)$ is singular. Therefore $f(x,y,z,t)=a_1x^3+a_2y^3+a_3z^3+a_4t^3+b_{34}z^2t+b_{43}t^2z$.
Unless $f_{A_3^{-1}}=f$, $V(f)$ is singular. Thus $f(x,y,z,t)=a_1x^3+a_2y^3+a_3z^3+a_4t^3$. Now $V(f)$ is nonsingular if and only if $a_1a_2a_3a_4\not=0$.

\section{$\SSS{5}$-invariant nonsingular cubic surfaces}
\setcounter{equation}{1}
A subgroup of $PGL_4(k)$ isomorphic to the symmetric group $\SSS{5}$ is one of three groups $C_{5!}$\Roman{equation},
\setcounter{equation}{2}
$C_{5!}$\Roman{equation},
\setcounter{equation}{3}
$C_{5!}$\Roman{equation} up to conjugacy \cite{mas}. We denote these groups by $G(1)$, $G(2)$ and $G(3)$, respectively. There exist
group isomorphisms $\varphi_i:\SSS{5}\rightarrow G(i)$ such that $(E_1)=\varphi_i((123))$, $(E_2)=\varphi_i((12)(34))$, $(E_3)=\varphi_i((12)(45))$
and $(F)=\varphi_i((12))$ generate $G(i)$ ($i\in [1,3]$) (cf \cite[(2.14) in chap.3]{suz}). Note that $(123)=(12)(23)$. Let $\omega=\frac{-1+i\sqrt{3}}{2}$. Then
$G(1)=\langle(E_1),(E_2),(E_3),(F)\rangle$, where
\begin{eqnarray*}
E_1&=&\left[\begin{array}{cccc}
            1&0&0&0\\
            0&1&0&0\\
            0&0&\omega&0\\
            0&0&0&\omega^2\end{array}\right],\
 E_2=\left[\begin{array}{cccc}
            1&0&0&0\\
            0&-\frac{1}{3}&\frac{2}{3}&\frac{2}{3}\\
            0&\frac{2}{3}&-\frac{1}{3}&\frac{2}{3}\\
            0&\frac{2}{3}&\frac{2}{3}&-\frac{1}{3}\end{array}\right],\
 E_3=\left[\begin{array}{cccc}
            \frac{1}{4}&-\frac{\sqrt{15}}{4}&0&0\\
            \frac{\sqrt{15}}{4}&\frac{1}{4}&0&0\\
            0&0&0&1\\
            0&0&1&0\end{array}\right],\\
 F&=&\frac{1+i}{\sqrt{2}}\left[\begin{array}{cccc}
                             1&0&0&0\\
                             0&1&0&0\\
                             0&0&0&1\\
                             0&0&1&0\end{array}\right].
\end{eqnarray*}
\begin{lemma} The $G(1)$-invariant nonsingular cubic surface is $V(f)$, where
\[
 f(x,y,z,t)=3\sqrt{15}x^3+10(y^3+z^3+t^3)-3\sqrt{15}xy^2-6(\sqrt{15}x+5y)zt.
\]
Moreover, $\Aut(V(f))=G(1)$.
\end{lemma}
\Proof
Let $F'=\frac{\sqrt{2}}{1+i}F$. Assume that $V(f)$ is a $G(1)$-invariant nonsingular cubic surface, where $f(x,y,z,t)$ is a cubic homogeneous
polynomial of the form as in the proof of Lemma 1.3. Since $E_1^3=E_4$, $f_{E_1^{-1}}=\omega^i f$ for some $i\in [0,2]$. If $i=1$, then
$f(x,y,z,t)=b_{13}x^2z+b_{23}y^2z+b_{34}z^2t+t^2(b_{41}x+b_{42}y)+c_4xyz$. Since $f_{F'^{-1}}\sim f$, we see $b_{13}=b_{23}=0$, hence
$V(f)$ is singular at $(1,0,0,0)$. Similarly $i\not=2$, namely $i=0$, so that $f(x,y,z,t)=a_1x^3+a_2y^3+a_3z^3+a_4t^3+b_1x^2y+b_2y^2x+c_1yzt+c_2xzt$.
We have $a_3a_4\not=0$, for $V(f)$ is nonsingular. If $a_1=a_2=0$, the condition $f_{E_3^{-1}}\sim f$ yields $b_1=b_2=0$, hence $V(f)$ is singular.
Thus $[a_1,a_2]\not=[0,0]$. Therefore $f_{F'^{-1}}=f$, hence $a_3=a_4$. Now $f_{E_3^{-1}}=f$, and $f_{E_3^{-1}}=-f$ does not hold.
The condition $f_{E_3^{-1}}=f$ is equivalent to $c_2=\sqrt{15}c_1/5$ and
\[
 \left[\begin{array}{cccc}
        -65&15\sqrt{15}&\sqrt{15}&-15\\
        15\sqrt{15}&-63&15&\sqrt{15}\\
        3\sqrt{15}&45&-93&13\sqrt{15}\\
        -45&3\sqrt{15}&13\sqrt{15}&-35\end{array}\right]
\left[\begin{array}{c}
      a_1\\
      a_2\\
      b_1\\
      b_2\end{array}\right]=0.
\]
The general $[a_1,a_2,b_1,b_2]$ satisfying the second condition above is $\lambda[3\sqrt{15},10,0,-3\sqrt{15}]$, for the rank of the above
matrix is equal to three. Consequently we may assume $f(x,y,z,t)=3\sqrt{15}x^3+10y^3+a_3(z^3+t^3)-3\sqrt{15}xy^2+c_1(yzt+\frac{\sqrt{15}}{5}xzt)$.
Now the condition $f_{E_3^{-1}}\sim f$ implies $f_{E_3^{-1}}=f$. Using the formula $(a+b+c)^3=a^3+b^3+c^3+3a^2(b+c)+3b^2(a+c)+3c^2(a+b)+6abc$,
we can write $f_{E_3^{-1}}(x,y,z,t)$ as
\begin{eqnarray*}
&&3\sqrt{15}x^3+(-\frac{10}{27}+\frac{16}{27}a_3-\frac{4}{27}c_1)y^3
+(\frac{80}{27}+\frac{7}{27}a_3-\frac{4}{27}c_1)z^3+(\frac{80}{27}+\frac{7}{27}a_3-\frac{4}{27}c_1)t^3\\
&&+y^2\{(-\frac{\sqrt{15}}{3}+\frac{4\sqrt{15}}{45}c_1)x+(\frac{60}{27}+\frac{12}{27}a_3+\frac{6}{27}c_1)(z+t)\}\\
&&+z^2\{(-\frac{4\sqrt{15}}{3}-\frac{2\sqrt{15}}{45}c_1)x+(-\frac{120}{27}+\frac{30}{27}a_3+\frac{6}{27}c_1)y+(\frac{240}{27}-\frac{6}{27}a_3+\frac{6}{27}c_1)t\}\\
&&+t^2\{(-\frac{4\sqrt{15}}{3}-\frac{2\sqrt{15}}{45}c_1)x+(-\frac{120}{27}+\frac{30}{27}a_3+\frac{6}{27}c_1)y+(\frac{240}{27}-\frac{6}{27}a_3+\frac{6}{27}c_1)z\}\\
&&+(-\frac{240}{27}-\frac{48}{27}a_3+\frac{3}{27}c_1)yzt+(-\frac{8\sqrt{15}}{3}+\frac{\sqrt{15}}{9}c_1)xzt\\
&&+(\frac{4\sqrt{15}}{3}+\frac{2\sqrt{15}}{45}c_1)xyt +(\frac{4\sqrt{15}}{3}+\frac{2\sqrt{15}}{45}c_1)xyz.
\end{eqnarray*}
Now $a_3=10$, $c=-30$ and $f_{E_3^{-1}}=f$ holds, as desired. We can show that $V(f)$ is nonsingular.\\

Let $\varepsilon=-\frac{\sqrt{5}-1}{4}+i\frac{\sqrt{10+2\sqrt{5}}}{4}$, hence $\ord(\varepsilon)=5$, $\alpha=\frac{-\sqrt{5}+1}{2}$,
$\beta=\alpha^2$, $\gamma=-\alpha$, and let
\[
 H=\left[\begin{array}{cccc}
          \varepsilon^4&0&0&0\\
          0&\varepsilon^2&0&0\\
          0&0&\varepsilon&0\\
          0&0&0&\varepsilon^3\end{array}\right],\
I=\left[\begin{array}{cccc}
          1&\alpha&\beta&\gamma\\
          \alpha&\beta&\gamma&1\\
          \beta&\gamma&1&\alpha\\
          \gamma&1&\alpha&\beta\end{array}\right].
\]
\begin{lemma} $(H)$ and $(I)$ generate in $PGL_4(k)$ a group $G(1)'$ conjugate to $G(1)$. The $G(1)'$-invariant nonsingular cubic
surface is $V(x^2y+y^2z+z^2t+t^2x)$. Moreover $\Aut(V(x^2y+y^2z+z^2t+t^2x))=G(1)'$.
\end{lemma}
\Proof
The transpositions $(i \;  i+1)$ ($i\in [1,4]$) generate $\SSS{5}$, while $(12345)(j \; j+1)(12345)^{-1}=(j+1\;  j+2)$ ($j\in [1,3]$). Therefore
$(12345)$ and $(12)$ generate $\SSS{5}$. Clearly $(123)(34)(45)=(12345)$. In view of the group isomorphism $\varphi_1:S_5\rightarrow G(1)$,
$\varphi_1((12345))=(E_1FE_2FE_3)$ and $\varphi_1((12))=(F)$. Let
\[
 K=E_1FE_2FE_3=\left[\begin{array}{cccc}
                      -\frac{1}{4}&\frac{\sqrt{15}}{4}&0&0\\
                      -\frac{\sqrt{15}}{12}&-\frac{1}{12}&\frac{2}{3}&\frac{2}{3}\\
                      \frac{\sqrt{15}\omega}{6}&\frac{\omega}{6}&\frac{2\omega}{3}&-\frac{\omega}{3}\\
                      \frac{\sqrt{15}\omega^2}{6}&\frac{\omega^2}{6}&-\frac{\omega^2}{3}&\frac{2\omega^2}{3}\end{array}\right],\
 T=\left[\begin{array}{cccc}
          1&0&0&0\\
          0&1&0&0\\
          0&0&1&-1\\
          0&0&1&1\end{array}\right].
\]
Now
\begin{eqnarray*}
&& T^{-1}FT=\left[\begin{array}{cccc}
                1&0&0&0\\
                0&1&0&0\\
                0&0&1&0\\
                0&0&0&-1\end{array}\right],\
T^{-1}KT=\left[\begin{array}{cccc}
               -\frac{1}{4}&\frac{\sqrt{15}}{4}&0&0\\
               -\frac{\sqrt{15}}{12}&-\frac{1}{12}&\frac{4}{3}&0\\
               -\frac{\sqrt{15}}{12}&-\frac{1}{12}&-\frac{1}{6}&-i\frac{\sqrt{3}}{2}\\
               -i\frac{\sqrt{5}}{4}&-i\frac{\sqrt{3}}{12}&-i\frac{\sqrt{3}}{6}&-\frac{1}{2}\end{array}\right], \\
&&\det(T^{-1}KT-\lambda E_4)=\lambda^4+\lambda^3+\lambda^2+\lambda+1=\Pi_{i=1}^{4}(\lambda-\varepsilon^i).
\end{eqnarray*}
If $\lambda$ is an eigenvalue of $T^{-1}KT$  and $x\in k^4$ satisfies $(T^{-1}KT-\lambda E_4)x=0$, then
\[
 [x_1,x_2,x_3,x_4]=x_1[1,\ (1+4\lambda)\frac{\sqrt{15}}{15},\ (1+\lambda+3\lambda^2)\frac{\sqrt{15}}{15},\ i(-\lambda^4+\lambda^3)\frac{\sqrt{3}}{5}].
\]
So $S^{-1}T^{-1}KTS=\diag[\varepsilon^4,\varepsilon^2,\varepsilon,\varepsilon^3]$, where
\[
 S=\left[\begin{array}{rrrr}
         1&1&1&1\\
         (1+4\varepsilon^4)\frac{\sqrt{15}}{15}&(1+4\varepsilon^2)\frac{\sqrt{15}}{15}&(1+4\varepsilon)\frac{\sqrt{15}}{15}&(1+4\varepsilon^3)\frac{\sqrt{15}}{15}\\
        (1+\varepsilon^4+3\varepsilon^3)\frac{\sqrt{15}}{15}&(1+\varepsilon^2+3\varepsilon^4)\frac{\sqrt{15}}{15}&(1+\varepsilon+3\varepsilon^2)\frac{\sqrt{15}}{15}&
  (1+\varepsilon^3+3\varepsilon)\frac{\sqrt{15}}{15}\\
        i(-\varepsilon+\varepsilon^2)\frac{\sqrt{3}}{5}&i(-\varepsilon^3+\varepsilon)\frac{\sqrt{3}}{5}&i(-\varepsilon^4+\varepsilon^3)\frac{\sqrt{3}}{5}&
  i(-\varepsilon^2+\varepsilon^4)\frac{\sqrt{3}}{5}\end{array}\right]
\]
so that $\det S=-i\frac{8\sqrt{3}}{25}\varepsilon^3(\varepsilon-1)^3(3\varepsilon^3+6\varepsilon^2+4\varepsilon+2)$.
We denote the $(i,j)$-cofactor of $S=[s_{ij}]$ by ${\tilde{s}}_{ij}$. Note that $S^{-1}\diag[1,1,1,-1]S=E_4+S^{-1}\diag[0,0,0,-2]S$ and that
the $i$-th row of $S^{-1}\diag[0,0,0,-2]S$ is equal to $-2{\tilde{s}}_{4i}/\det S\times$(4-th  row of $S$).
Denote $TS\diag[1,\varepsilon^2,\varepsilon^3,\varepsilon]$ by $S'$. By computation we get
the components of $(3\varepsilon^3+6\varepsilon^2+4\varepsilon+2)S'^{-1}FTS$' as
follows.
\[
\left[\begin{array}{cccc}
      2\varepsilon^3+4\varepsilon^2+3\varepsilon+1& 2\varepsilon^4+\varepsilon^3+2& -\varepsilon^4+\varepsilon^2& \varepsilon^3+2\varepsilon^2+2\varepsilon\\
    2\varepsilon^4+\varepsilon^3+2& -\varepsilon^4+\varepsilon^2& \varepsilon^3+2\varepsilon^2+2\varepsilon& 2\varepsilon^3+4\varepsilon^2+3\varepsilon+1\\
    -\varepsilon^4+\varepsilon^2& \varepsilon^3+2\varepsilon^2+2\varepsilon& 2\varepsilon^3+4\varepsilon^2+3\varepsilon+1& 2\varepsilon^4+\varepsilon^3+2\\
    \varepsilon^3+2\varepsilon^2+2\varepsilon& 2\varepsilon^3+4\varepsilon^2+3\varepsilon+1& 2\varepsilon^4+\varepsilon^3+2& -\varepsilon^4+\varepsilon^2
     \end{array}\right].
\]
Moreover, since $(2\varepsilon^3+4\varepsilon^2+3\varepsilon+1)^{-1}=(7\varepsilon^3+4\varepsilon^2+\varepsilon+8)/5$, we can verify
\begin{eqnarray*}
&&(2\varepsilon^3+4\varepsilon^2+3\varepsilon+1)^{-1}(2\varepsilon^4+\varepsilon^3+2)=\varepsilon^3+\varepsilon^2+1=\alpha,\\
&&(2\varepsilon^3+4\varepsilon^2+3\varepsilon+1)^{-1}(-\varepsilon^4+\varepsilon^2)=\varepsilon^3+\varepsilon^2+2=\alpha^2,\\
&&(2\varepsilon^3+4\varepsilon^2+3\varepsilon+1)^{-1}(\varepsilon^3+2\varepsilon^2+2\varepsilon)=-\varepsilon^3-\varepsilon^2-1=-\alpha.
\end{eqnarray*}
 We have shown that $S'^{-1}KS'=H$ and $(S'^{-1}FS')=(I)$.

Finally we shall show that a cubic homogeneous polynomial $f(x,y,z,t)$ of the form as in the proof of Lemma 1.3 such that $f_{H^{-1}}\sim f$,
$f_{I^{-1}}\sim f$ and that $V(f)$ is nonsingular is proportional to $x^2y+y^2z+z^2t+t^2x$. Since $H^5=E_4$, $f_{H}^{-1}=\varepsilon^i f$ for
some $i\in [0,4]$. If $i=1$, then
$f(x,y,z,t)=a_2y^3+b_{14}x^2t+b_31z^3x+c_1yzt$, hence $V(f)$ is singular at $(0,0,0,1)$. Similarly, unless $i=0$, $V(f)$ has a singular point.
If $i=0$, then $f(x,y,z,t)=b_1x^2y+b_2y^2z+b_3z^2t+b_4t^2x$. Denoting the coefficients of $x^3$, $y^3$, $z^3$ and $t^3$ in $f_{I^{-1}}$ by

$b_1'$, $b_2'$, $b_3'$ and $b_4'$, respectively, and noting $\alpha^2=\alpha+1$, we obtain
\[
 b'=\alpha \left[\begin{array}{cccc}
                 1&\alpha^3&-\alpha^4&\alpha\\
                 \alpha^3&-\alpha^4&\alpha& 1\\
                 -\alpha^4&\alpha&1&\alpha^3\\
                 \alpha&1&\alpha^3&-\alpha^4\end{array}\right]b.
\]
Since $b'=0$ and the rank of the matrix involved is equal to three, it follows that $b_1=b_2=b_3=b_4\not=0$. Now let $f=x^2y+y^2z+z^2t+t^2x$.
In order to see $f_{I^{-1}}=5(4\alpha+3)f$, we note that denoting $(x+\alpha y+\beta z+\gamma t)^2(x+\alpha x+\beta y+\gamma z+t)$ by $g(x,y,z,t)$,
we have  $f_{I^{-1}}(x,y,z,t)=g(x,y,z,t)+g(t,x,y,z)+g(z,t,x,y)+g(y,z,t,x)$. \\

The group $G(2)$  is generated by $(E_1)$, $(E_2)$, $(E_3)$ and $(F)$, where
\begin{eqnarray*}
E_1&=&\left[\begin{array}{cccc}
            1&0&0&0\\
            0&1&0&0\\
            0&0&\omega&0\\
            0&0&0&\omega^2\end{array}\right],\
 E_2=\frac{1}{\sqrt{3}}\left[\begin{array}{cccc}
            1&0&0&\sqrt{2}\\
            0&-1&\sqrt{2}&0\\
            0&\sqrt{2}&1&0\\
            \sqrt{2}&0&0&-1\end{array}\right],\
 E_3=\left[\begin{array}{cccc}
            \frac{\sqrt{3}}{2}&\frac{1}{2}&0&0\\
            \frac{1}{2}&-\frac{\sqrt{3}}{2}&0&0\\
            0&0&0&1\\
            0&0&1&0\end{array}\right],\\
F&=&\left[\begin{array}{cccc}
                             0&1&0&0\\
                             -1&0&0&0\\
                             0&0&0&1\\
                             0&0&-1&0\end{array}\right].
\end{eqnarray*}
\begin{lemma}
There exist no $G(2)$-invariant nonsingular cubic surfaces.
\end{lemma}
\Proof
Let $f(x,y,z,t)$ be a $G(2)$-invariant cubic homogeneous polynomial of the form as in the proof of Lemma 1.3. Since $F^4=E_4$, $f_{F^{-1}}=i^j f$
for some $j\in[0,3]$. If $j=0$ or $j=2$, then $f=0$. If $j=1$, then $f(x,y,z,t)$ has the form
\begin{eqnarray*}
&&a_1(x^3-iy^3)+a_3(z^3-it^4)\\
&&+x^2(b_{12}y+b_{13}z+b_{14}t)+y^2(ib_{12}x+ib_{14}z-ib_{13}t)\\
&&+z^2(b_{31}x+b_{32}y+b_{34}t)+t^2(ib_{32}x-ib_{31}y+ib_{34}z)\\
&&+c_1yzt-ic_1xzt+c_3xyt-ic_3xyz.
\end{eqnarray*}
Since $E_3^2=E_4$, $f_{E_3^{-1}}=\pm f$.  The coefficients of $t^3$, $t^2z$, $t^2x$ and $t^2y$ of $f_{E_3^{-1}}$ (resp. $\pm f$) are
\[
 a_3,\ b_{34},\ \frac{\sqrt{3}}{2}b_{31}+\frac{1}{2}b_{32},\ \frac{1}{2}b_{31}-\frac{\sqrt{3}}{2}b_{32}\
({\rm resp.\ }\mp ia_3,\ \pm ib_{34},\ \pm ib_{32},\ \mp ib_{31}),
\]
hence $a_3=b_{34}=b_{31}=b_{32}=0$. Consequently $V(f)$ is singular at $(0,0,0,1)$.\\

Let $\nu_1=\sqrt{\frac{3}{8}}+i\sqrt{\frac{5}{8}}$ and $\nu_2=\sqrt{\frac{3}{8}}-i\sqrt{\frac{5}{8}}$
The group $G(3)$ is generated by $(E_1)$, $(E_2)$, $(E_3)$ and $(F)$, where
\begin{eqnarray*}
E_1&=&\left[\begin{array}{cccc}
            \omega&0&0&0\\
            0&\omega&0&0\\
            0&0&\omega^2&0\\
            0&0&0&\omega^2\end{array}\right],\
 E_2=\frac{1}{\sqrt{3}}\left[\begin{array}{cccc}
            1&0&0&\sqrt{2}\\
            0&1&\sqrt{2}&0\\
            0&\sqrt{2}&-1&0\\
            \sqrt{2}&0&0&-1\end{array}\right],\
 E_3=\left[\begin{array}{cccc}
            0&0&0&\nu_1\\

            0&0&\nu_2&0\\

            0&\nu_1&0&0\\

            \nu_2&0&0&0\end{array}\right],\\
 F&=&\left[\begin{array}{cccc}

                             0&0&0&1\\
                             0&0&-1&0\\
                             0&1&0&1\\
                             -1&0&0&0\end{array}\right].
\end{eqnarray*}

\begin{lemma}
There exist no $G(3)$-invariant nonsingular cubic surfaces. More precisely, there exist no $G(3)$-invariant cubic forms except zero.
\end{lemma}
\Proof
Let $E_1'=\omega^{-1}E_1$ and $f(x,y,z,t)$ a $G(3)$-invariant cubic form. If $f_{E_1'^{-1}}=\omega f$, then $f$ has the form
\[
 x^2(b_{13}z+b_{14}t)+y^2(b_{23}z+b_{24}t)+c_3xyt+c_4xyz,
\]
hence the condition $f_{F^{-1}}\sim f$ implies $f=0$. Similarly, if $f_{E_1'^{-1}}=\omega^2 f$, then $f=0$. Now assume $f_{E_1'^{-1}}= f$, so that
$f(x,y,z,t)=a_1x^3+a_2y^3+a_3z^3+a_4t^3+b_1x^2y+b_2y^2x+b_3z^2t+b_4t^2z$. Since $f_{F^{-1}}=i^j f$ for some $j\in [0,3]$, we obtain
\begin{eqnarray*}
\left[\begin{array}{cc}
       i^j&1\\
       -1&i^j\end{array}\right]\left[\begin{array}{c}
                                      a_1\\
                                      a_4\end{array}\right]=0,\
           \left[\begin{array}{cc}
                 i^j& -1\\
                 1&i^j\end{array}\right]\left[\begin{array}{c}
                                              a_2\\
                                              a_3\end{array}\right]=0,
                                               \left[\begin{array}{cc}
                                                           i^j&-1\\
                                                           1&i^j\end{array}\right]\left[\begin{array}{c}
                                                                                          b_1\\
                                                                                          b_4\end{array}\right]=0,\
                                                                                                \left[\begin{array}{cc}
                                                                                                      i^j& 1\\
                                                                                                      -1&i^j\end{array}\right]\left[\begin{array}{c}
                                                                                                                                              b_2\\
                                                                                                                                            b_3\end{array}\right]=0.
\end{eqnarray*}
Therefore, unless $j=1$ or $j=3$, $f(x,y,z,t)=0$. Even if $j=1$ or $j=3$, the condition $f_{E_3^{-1}}\sim f$ implies $f=0$, because
$\{i,-i\}\cap\{\nu_1,\nu_2,\nu_1^3,\nu_2^3\}=\emptyset$.\\

\begin{remark} Let $G=\Aut(V(x^2y+y^2z+z^2t+t^2x))$, $\Sigma$ the field automorphism of $\QQ(\varepsilon)$ such that $\varepsilon^\Sigma=\varepsilon^2$.\\
$(1)$ Denoting by $e_i$ be the $i$-th column of the unit matrix $E_4$, let $J=[e_2,e_3,e_4,e_1]$. Then $(J)\in G=G(1)'$.\\
$(2)$ Define $I'\in GL_4(k)$ as $I$, using $\alpha'=\frac{\sqrt{5}+1}{2}$, $\beta'={\alpha'}^2$ and $\gamma'=-\alpha'$ in place of
$\alpha$, $\beta$ and $\gamma$, respectively. Then $H^\Sigma=H^2$ and $I^\Sigma=I'$, for ${\sqrt{5}}^\Sigma=-\sqrt{5}$. On the other hand
$JHJ^{-1}=H^2$ and $JIJ^{-1}=\beta I'$. Consequently there exists uniquely a group isomorphism $\psi:\SSS{5}\rightarrow PGL_4(k)$ such that
$\psi((12345))=(H^2)$ and $\psi((12))=(I')$. The projective representations of $\SSS{5}$ $\psi$ and $\varphi_1$ are equivalent.
\end{remark}

%
\end{document}